\documentclass[twoside,leqno,fleqn]{article}
\usepackage{amsmath}
\usepackage{amssymb}
\usepackage[shortlabels]{enumitem}
\usepackage[onehalfspacing]{setspace}

\newcommand{\AUTOR}{C. Ga\ss ner}
\newcommand{\TITEL}{AC and the Independence of the Law of Trichotomy in Second-Order Henkin Logic}
\def\S{\Sigma}
\def\I{{\cal I}}
\def\G{\mathfrak{G}}

\newcommand{\bbbn}{\mathbb{N}} 
\newcommand{\mbm}[1]{\mbox{\boldmath{$#1$}}}
\newcommand{\mbmss}[1]{\mbox{\scriptsize\boldmath{$#1$}}}
\newcommand{\qed}{\hfill{$\Box$}} 

\newtheorem{satz}{Satz}[section] 
\newtheorem{lemma}[satz]{Lemma}
\newtheorem{proposition}[satz]{Proposition}
\newtheorem{corollary}[satz]{Corollary}
\newtheorem{theorem}[satz]{Theorem}
\newtheorem{remark}[satz]{Remark}
\newtheorem{idea}[satz]{Idea}
\newtheorem{question}[satz]{Question}
\counterwithout{equation}{section}

\begin{document}
\pagestyle{headings} 
\thispagestyle{empty}
\begin{center} {\Large\bf AC and the Independence of the Law of Trichotomy in Second-Order Henkin Logic \footnote{Following on from a presentation at the Colloquium Logicum 2022, the first part of this paper deals with a proof whose implications were discussed in Leeds and whose ideas were presented in detail at a meeting in Kloster 2024. I would like to thank the staff of the International Office of the University of Greifswald for their support. My thanks go also to the organizers of the Colloquium Logicum 2022 in Konstanz.
}}\vspace{0.6cm}\\{\bf Christine Ga\ss ner} {\vspace{0.2cm}\\University of Greifswald, Germany, 2024\\ gassnerc@uni-greifswald.de}\\\end{center} 

\begin{abstract} 
This paper focuses on the set HAC of 1-1 Ackermann axioms of choice in second-order predicate logic with Henkin interpretation (HPL). To answer a question posed by Michael \text{Rathjen}, we restrict the proof that the basic Fraenkel model of second order is a model of all n-m Ackermann axioms to the case where the Ackermann axioms are in HAC. In the second part, we show the independence of Hartogs' version of the law of trichotomy (TR) from HAC in HPL. A generalization of the latter proof implies the independence of TR from all Ackermann axioms in HPL. We conclude the paper with an open problem.
\end{abstract}

\section{Introduction}\label{SectionEinf}

Friedrich M.\,\,Hartogs presented a proof in \cite{Hartogs} that implies the equivalence of the law of trichotomy
 and the well-ordering principle in the usual set theory. In \cite[p.\,128]{1c}, G\"unter Asser pointed out that, for second-order formulations of principles such as the Axiom of Choice (AC), Zorn's lemma (ZL), the law of trichotomy (TR), and the well-ordering principle (WO), many of the classical proofs cannot be reconstructed in second-order predicate logic with Henkin interpretation (HPL). This prompted us to examine the known proofs in more detail and to look for counterexamples. On the theoretical basis of \cite{1c}, the strength of the principles of choice and the relationships between these and other second-order formulas in HPL was analyzed in \cite{Gass84}. Our paper \cite{Gass94} and the tables in \cite[p.\,17--18]{Gass24D} contain comprehensive overviews that are useful to classify these statements. We therefore also know that the second-order versions $AC^{1,1}$ and $TR^1$ are HPL-independent of each other and that $WO^1$ is HPL-independent from $AC^{1,1}$ and from $TR^1$. This approach was continued in \cite{Gass24B,Gass24C}. Since the Ackermann axioms in HAC are stronger than $AC^{1,1}$ in HPL (see \cite[Theorem 3.2]{Gass24C}), we will also address a question posed by Stephen Mackereth in \cite{Mackereth} and show that $TR^1$ is HPL-independent from the Ackermann axioms. The proofs show that the 1-1 Ackermann axioms belonging to the set HAC hold in several permutation models. This also answers a question posed by Michael \text{Rathjen}. 

Here, we consider common properties of the permutation models $\S_0$ in $ {\sf struc}_{\rm pred}^{({\rm m})}( \bbbn) $ and $k\S_0$ in $ {\sf struc}_{\rm pred}^{({\rm m})}(\{1,\ldots,k\}\times \bbbn)$ introduced in \cite{Gass84} and \cite{Gass24D} (for $k=2$), respectively, and defined by the Fraenkel-Mostowski-Specker-Asser method (that was discussed in \cite{1c,Gass84,Gass94,Gass24D}). The components used in the construction of these models are an infinite set of individuals or a finite number of pairwise disjoint copies of this set, sets of $n$-ary predicates ($n\geq1$) that can be defined by using finite sets of individuals and the identity and by applying set-theoretical operations, groups of permutations, systems of finite sets of individuals, and the like. Consequently, it is possible to use a metatheory based on just a few axioms. The axiom of choice is not necessary for the construction of $\S_0$ and $k\S_0$. 

To give a simplified proof of the HPL-independence of $WO^n$ from the 1-1 Ackermann axioms that is derived from the proof of Proposition 4.14 in \cite{Gass24B}, we will use notations similar to those summarized in \cite[Table 1]{Gass24B}. To make it easier to compare the proofs of $\S_0\models choice_h^{1,1}$ and $\S_0\models choice_h^{n,m}$ and to show the parallels between the proof presented here and the proof of the more general statement in \cite[Prop. 4.14]{Gass24B}, we will include an overview of the corresponding notations used here in Table \ref{overview_nota_1}. Moreover, for the evaluation of the model $k\S_0$, some notations will be modified again as indicated in Table \ref{overview_nota_2}.

When formalizing the Ackermann axioms of choice and other statements, we will also use the following variables. For given $n\geq 1$ and $m\geq 1$, $\mbm{x}$ stands for the $n$ variables $x_1,\ldots,x_n$, $\mbm{y}$ stands for the $m$ variables $y_1,\ldots,y_m$, and so on. For $m=1$ and $n=1$, we use $x$ instead of $\mbm{x}$ and $y$ instead of $\mbm{y}$. Let $A$ and $B,B_1, B_2, \ldots$ be variables of sort $n$, let $D$ be of sort $m$, and let $R, R_1,\ldots $ and $S, S_0, \ldots $ be of sort $n+m$. The set ${\cal L}^{(2)}_{\mbmss{x},D}$ contains all second-order formulas in which the individual variables given by $\mbm{x}$ and $D$ occur only free, ${\cal L}^{(2)}$ contains all second-order formulas, and so on. The string $=_{\rm df}$ used sometimes in equations indicates that such an equation is also a definition of the term on the left-hand side which is the centre of attention at that moment. We will consider the following axioms ($n,m\geq 1$). 

\vspace{0.2cm}
\noindent{\bf The set $choice_h^{(2)}$ of all Ackermann axioms in HPL} 

\nopagebreak 
\noindent\fbox{\parbox{11.7cm}{
\vspace*{0.1cm}
$\begin{array}{ll}choice_h^{(2)}&=_{\rm df}\bigcup_{n,m\geq 1 }choice_h^{n,m}\\
choice_h^{n,m}&=_{\rm df}\{choice_h^{n,m}(H)\mid H \mbox{ is a formula in ${\cal L}^{(2)}_{\mbmss{x},D}$} \}\\
 choice_h^{n,m}(H)\!\!&=_{\rm df} \! \forall \mbm{x} \exists D H(\mbm{x},D) \!\to \exists S \forall \mbm{x } \exists D (\forall\mbm{y}(D \mbm{y} \!\leftrightarrow S\mbm{x}\mbm{y}) \land H(\mbm{x}, D )) 
\end{array}$
} } 

\vspace{0.2cm}
\noindent{\bf The set {\rm HAC} of Ackermann axioms} 

\nopagebreak 
\noindent\fbox{\parbox{11.7cm}{
\vspace*{0.1cm}
$\begin{array}{ll}
{\rm HAC}\hspace{0.14cm}&=_{\rm df}\quad \{choice_h^{1,1}(H)\mid H \mbox{ is a formula in ${\cal L}^{(2)}_{x,D}$} \}\\
\end{array}$
} } 

\vspace{0.2cm}
\noindent Thus, {\rm HAC} is the set of the 1-1 Ackermann axioms. Let us remark that the independence of $WO^n$ and HAC from each other was also investigated by Benjamin Siskind, Paolo Mancosu, and Stewart Shapiro in \cite{Siskind}.

Here, we will also consider second-order statements $TR^n$ for $n$-ary predicates ($n\geq 1$), whose first-order counterpart $TR$ is equivalent to the axiom $AC$ of choice in ZF. 

\vspace{0.2cm}
\noindent{\bf The Hartogs axioms of trichotomy for $n$-ary predicates}

\nopagebreak 
\noindent\fbox{\parbox{11.7cm}{

\vspace*{0.1cm}
$TR^n\hspace{0.35cm}=_{\rm df} \quad \forall A \forall B ( [A\lesssim B] \lor [B\lesssim A])$.
}}

\vspace{0.2cm}

\noindent $TR^n$ is a formulation of the law of trichotomy, $TR$, in second-order predicate logic. $[A\lesssim B]$ stands for the statement that there is a one-to-one function of $A$ into $B$. $TR^n$ states that for any two $n$-ary predicates there is a one-to-one function --- given in the form of a $2n$-ary predicate --- of one of these predicates into the other. According to \cite[p.\,47, Section V, Satz 1.2]{Gass84} and \cite{Gass94}, $\S_0$ is a model of $TR^1$. Thus, $TR^1\to WO^1$ does not hold in HPL. Because of $^hax^{(2)}\vdash WO^1\to TR^1$, the independence of $TR^1$ from HAC in HPL would imply the HPL-independence of $WO^1$ from HAC, but the HPL-independence of $WO^1$ does not imply the HPL-independence of $TR^1$. To show that the second-order formulas $TR^n$ ($n\geq 1$) are independent from HAC, we will use the permutation model $k\S_0$ ($k\geq 2$). Consequently, this will once more confirm the HPL-independence of $WO^1$ from HAC. Just as the first proof in this article deals with the axioms in HAC instead of all Ackermann axioms considered in the original proof, the second proof can be generalized in such a way that it is possible to replace the axioms in HAC by any $n$-$m$ Ackermann axioms. 

For more details and further references, see \cite{Gass24A}, \cite{Siskind}, and \cite{Mackereth}. We would like to point out that our research was also influenced by \cite{Ack110,Ack111, Felgner, FA22b, FA22a, FA37, Henk, Jech73, Jech03, Kame, M39, 13, Specker, Vaeae}. 

In Section \ref{Section3}, we prove that the basic Fraenkel model II is a model of the axioms belonging to {\rm HAC}. In Section, \ref{Section4} we consider further permutation models in which the axioms belonging to ${\rm HAC}$ hold. The models are defined in Section \ref{Section2_1}. The proof ideas are described in Sections \ref{Sect_Idea1} and \ref{Sect_Idea2}. Roughly speaking, the main idea of the proofs of $\S_0\models choice_h^{1,1}(H)$ and $k\S_0\models choice_h^{1,1}(H)$ for all formulas $H\in {\cal L}^{(2)}_{x,D}$ is to show the existence of the required predicates by decomposing the problem into a finite number of subproblems defined by a finite stabilizer for $H$. This allows to define a suitable choice function derived from a choice function on a finite family of sets. Some of the necessary prerequisites can be found in Sections \ref{Sect_FiniteSupports} and \ref{Sect_FiniteChoice}. Adequate decompositions and the consequences of swapping two individuals by permutations are discussed in Sections \ref{Decomp} and \ref{DefinPerm}. Sections \ref{sectionConseq} and \ref{section5} contain conclusions and an open problem.

\section{Finite, infinite, basic, and standard structures} \label{Section2}

In our metatheory, we associate the term {\sf finite set} with sets such as the following: $\{0\}$, $\{0,1\}$, $\{0,1,2\}$, \ldots, and the like. We will use the axiom of infinity which guarantees that there is an infinite set. However, for simplicity matters, we assume that there is a set that contains at least all elements of the sets $\{0\}$, $\{0,1\}$, $\{0,1,2\}$, \ldots. A smallest set of this kind should be the {\sf infinite set} which we will use. We denote it by $\bbbn$ and assume $0<1<2<\cdots$ and $0\leq 1\leq 2\leq \cdots$. Moreover, we will use the principle of induction. Any other set is {\sf finite} if there is a 1-1 function of this set onto one of the sets $\{0\}$, $\{0,1\}$, $\{0,1,2\}$, \ldots. For simplicity matters, let the (first) element of the first non-empty intersection of a co-finite subset $\bbbn\setminus P$ and one of the sets $\{0\}$, $\{0,1\}$, $\{0,1,2\}$, \ldots (in this order) be called the {\sf minimum} of the set $\bbbn\setminus P$ and denoted by $\min(\bbbn\setminus P)$. The symbol $\in$ is used to say that an object {\sf belongs to} a domain or a set. The symbol $\Rightarrow$ stands for {\sf implies} ({\sf that}).

We are interested in an infinite domain or set of individuals, finite stabilizers, finite families, finite structures, finite unions of structures, and so on.

\subsection{The models used}\label{Section2_1}

Let $I$ be any non-empty set, $\G_{\sf 1}^I$ be the group $perm(I)$ containing all permutations of $I$, $\G$ be a subgroup of $\G_{\sf 1}^I$, and $\I_{\sf 0}^I$ be the set of all finite subsets of $I$. For each predicate $\alpha \in {\rm pred}(I)$, let ${\rm sym}_{\G}(\alpha)$ be defined by ${\rm sym}_{\G}(\alpha)=\{ \pi \in {\G}\mid \alpha^\pi=\alpha\}$. For any $P\subseteq I$, let $\G(P)=\{\pi \in {\G}\mid \forall \xi (\xi \in P\Rightarrow \pi(\xi)=\xi) \}$. For any $n\geq1$, let $J_{n}(I,\G,\I_{\sf 0}^I)=\{\alpha \in {\rm pred}_n (I)\mid (\exists P\in \I_{\sf 0}^I)({\rm sym}_{\G}(\alpha)\supseteq \G(P))\}$ and let $\S(I,\G,\I_{\sf 0}^I)$ be the predicate structure $ (J_{n})_{n\geq 0}$ with the domains $J_{0}=I$ and $J_{n}=J_{n}(I,\G,\I_{\sf 0}^I)$ for $n\geq 1$. We will consider the basic Fraenkel model $\S_0$ of second order given by $\S_0= \S(\bbbn,\G_{\sf 1}^\bbbn,\I_{\sf 0}^I)$, the permutation model $k\S_0$
given by 

$
\begin{array}{ll}
k&\geq 1,\\
I&=\{1,\ldots,k\}\times \bbbn,\\ 
I_j&=\{j\}\times \bbbn$ \hspace{5.1cm} for all $j\in\{1,\ldots,k\},\\
\G&=\{\pi \in \G_{\sf 1}^I\mid (\forall j\in \{1,\ldots,k\} ) (\forall \xi \in I_j) (\pi(\xi)\in I_j)\},\\
k\S_0&=\S(I, \G,\I_0^I),
\end{array}$

\noindent and the standard structures $\S_1^{\bbbn}$ and $\S_1^{\{1,\ldots,k\}}$ ($k\geq 1$) given by $\S_1^I= \S(I,\G_{\sf 0}^I,\I_{\sf 0}^I)$ where $\G_{\sf 0}^I$ is the subgroup of $\G_{\sf 1}^I$ containing only the identity permutation of $I$. Note that $\S_1^{\{1,\ldots,k\}}$ is a substructure of both Henkin-Asser structures $\S_0$ and $\S_1^{\bbbn}$ and it is the only Henkin-Asser structure in ${\sf struc}_{\rm pred}^{({\rm m})}(\{1,\ldots,k\})$.

\subsection{Finite stabilizers}\label{Sect_FiniteSupports}

Based on the theoretical foundations provided by Günter Asser in \cite{1c}, we consider stabilizers for formulas. The following lemma follows from \cite[Prop. 3.7]{Gass24D} (and \cite[Prop. 2.21]{Gass24A}). \cite{Gass24D} includes also a proof of \cite[Prop. 3.7]{Gass24D}. 

\begin{lemma}[Finite supports and stabilizers for formulas]\label{EndlichesPfuerH} Let $H( x ,D)$ be in ${\cal L}^{(2)}_{ x ,D}$. Let $\S$ be any structure of the form $ \S(I,{\G},\I_{\sf 0}^I)$ for any subgroup $\G\subseteq {\G}_{\sf 1}^I$ and let $f$ be in ${\rm assgn}(\S)$. Then, there is a finite set $P$ in $\I_{\sf 0}^I$ such that \begin{equation}\label{Stabil}\S_{f\langle{x \atop \xi}{D\atop\delta}\rangle}(H( x ,D))=\S_{f\langle{ x \atop \pi(\xi)}{D\atop\delta^\pi}\rangle}(H( x ,D))\end{equation} holds for all $\xi\in I$, for all predicates $\delta\in J_1(I,{\G},\I_{\sf 0}^I)$, and for all $\pi\in \G(P)$.
\end{lemma}
\noindent {\bf Proof.}
Let $r,s\geq 0$. If $s>0$, then let $x_{i_1},\ldots,x_{i_s}$ be the individual variables that may additionally occur free in $H(x,D)$ such that $x\not\in \{x_{i_1},\ldots,x_{i_s}\}$ holds. If $r>0$, then let $A_{j_1}^{n_1}, \ldots, A_{j_r}^{n_r}$ be the predicate variables that may additionally occur free in $H(x,D)$ such that $D\not\in \{A_{j_1}^{n_1}, \ldots, A_{j_r}^{n_r}\}$ holds. Assume that $H(x,D)$ does not contain further free variables. Let $P'=\{f(x_{i_1}) , \ldots, f(x_{i_s})\}$ if $s>0$ and let $P'=\emptyset$ if $s=0$. If $r=0$, then let $P=P'$. If $r>0$, then, for every $i\in\{1,\ldots,r\}$, let $P_i$ be a finite subset of $I$ such that ${\rm sym}_{\G}(\alpha_i)\supseteq \G(P_i)$ holds for $\alpha_i=f(A_{j_i}^{n_i})$. Moreover, let $P=P'\cup P_1\cup \cdots \cup P_r$. Thus, in any case, we have the equation (\ref{Stabil}) for each $\pi \in \G(P)$. \qed

\subsection{Choice functions and finite families}\label{Sect_FiniteChoice}

From the point of view that we want to analyze a proof (of \cite[Prop. 4.14]{Gass24B}) and confirm that the Ackermann axiom $choice^{1,1}_h(H)$ holds in $\S_0$ for any second-order formula $H( x ,D)$ belonging to ${\cal L} ^{(2)}_{ x ,D}$,
the properties \ref{Antec} and \ref{Conseq} considered in the next lemma are important for us.
 \ref{Antec} and \ref{Conseq} reflect the meaning of the antecedent and the consequent, respectively, of the Ackermann axioms in $choice^{1,1}_h$. Recall that $\widetilde{\alpha}$ is the set $ \{\xi \mid \alpha(\xi)=true\}$ and that we sometimes write $\alpha$ for $\widetilde{\alpha}$, and so on. 

\begin{lemma}[The assumption $\forall x \exists D H( x ,D)$ and choice functions]\label{Zhg_CG_sigma}\hfill Let \newline $\S$ be any predicate structure $(J_n)_{n\geq 0}$ in ${\sf struc}_{\rm pred}^{({\rm m})}(J_0)$ and $f$ be any assignment in $\S$. For any formula $H$ in ${\cal L} ^{(2)}_{ x ,D}$, let 
\[ {\cal A}^H=\{( \xi , \delta )\in J_0 \times J_1\mid \S\models_{f\langle{x \atop \xi}{D\atop\delta}\rangle} H\}\]
 and, for any $\xi \in J_0 $, let
\[{\cal A}^H_{\xi}=\{ \delta \in J_1\mid ( \xi , \delta )\in {\cal A}^H\}.\] 

\begin{enumerate}[label={\rm(\arabic*)}]
 \labelwidth0.7cm \leftmargin0cm \itemsep2pt plus1pt
\topsep1pt plus1pt minus1pt
\labelsep4pt \parsep0.5pt plus0.1pt minus0.1pt
\item\label{Antec} $\S\models_f\forall x \exists D H( x ,D)$ implies that, for all individuals $\xi $ in $J_0$, there is a non-empty uniquely determined subset ${\cal A}^H_{\xi}\subseteq J_1$.
\item \label{Conseq} If there is a choice function $\varphi$ on the family $({\cal A}^H_{\xi})_{ \xi \in J_0 }$ such that $\varphi ( {\cal A}^H_{\xi})\in {\cal A}^H_{\xi}$ holds for all $\xi \in J_0 $ and if, moreover, the $(1+1)$-ary predicate $\sigma_{\varphi}$ defined by $\widetilde{\sigma_{\varphi}}=\{(\xi , \eta)\mid \xi \in J_0 \,\,\&\,\, \eta \in\varphi ( {\cal A}^H_{\xi})\}$ is in $J_{1+1}$, then we have $\S\models_{f\langle{S\atop\sigma_{\varphi}} \rangle}\forall x \exists D ( \forall y (D y \leftrightarrow S x y ) \land H( x , D ))$. 
\end{enumerate}
\end{lemma}

Besides this observation, we have the following for non-empty finite families $\subseteq J_0\times {\cal P}(J_1)$ (and the corresponding statements for all finite families $\subseteq J_0^n\times {\cal P}(J_m)$ with $n,m\geq 1$) since we know that the existence of choice functions $\phi_{s-1}$ ($s\geq 1$) on finite systems of sets can be proved by induction on the number $s$ of sets in such systems as explained by Abraham A.\,\,Fraenkel in \cite{FA37}. Recall that ${\cal P}(J_m)$ denotes the power set of $J_m$. 

\begin{lemma}[$choice_h^{1,1}(H)$ restricted to finite families $({\cal A}^H_{\xi})_{ \xi \in \alpha}$]\label{endlSigma} \hfill For any \linebreak Henkin-Asser structure $\S=(J_n)_{n\geq 0}$, any finite $\alpha \in J_1$, any $f$ in ${\rm assgn}(\S)$, and $H( x ,D)\in {\cal L} ^{(2)}_{ x ,D}$, there holds
\[\S\models_{f\langle {A\atop \alpha}\rangle}\forall x \exists D (A x \to H( x ,D))\to \exists S \forall x \exists D ( A x \to \forall y (D y \leftrightarrow S x y ) \land H( x , D )).\] 
\end{lemma}

This means that the finite Henkin structure $(J_n)_{n\geq 0}$ of the form $\S_0^{\{1,\ldots,k\}}$ with $k\geq 1$ is a model of $choice_h^{1,1}(H)$. This result can be generalized since, for any $n,m\geq 1$, any relevant family $({\cal A}^H_{\mbmss{\xi}})_{ \mbmss{\xi} \in J_0^n}$ that comes into question is included in $\widetilde{\alpha_0} \times {\cal P}(J_m)$ for the finite predicate $\alpha_0=_{\rm df}\{1,\ldots,k\}^n$. On the other hand, $1<2<\cdots <k$ induces a well-ordering on $\widetilde{\alpha_0}$. Consequently, we have the following statements.

\begin{proposition}[$choice_h^{(2)}\!$ and $WO^n$ hold in finite Henkin structures]\label{IndepWO} \hfill

\noindent For any $k \geq 1$ and $n\geq 1$, we have 

 $\S_0^{\{1,\ldots,k\}}\models choice_h^{(2)}$,

 $\S_0^{\{1,\ldots,k\}} \models WO^{n}$. 
\end{proposition}

\vspace{0.1cm}

 \begin{remark}[The ${\rm HPL}$-independence of $WO^n$]\label{WOinFinStr}
Proposition \ref{IndepWO} confirms \linebreak that $WO^n$ is {\rm HPL}-independent from $ choice_h^{(2)}$ because $\S_0\models choice_h^{(2)}\cup\{\neg WO^n\}$ was proved in \cite{Gass24B} and will be confirmed below for ${\rm HAC}\cup\{\neg WO^n\}$.
\end{remark}

\vspace{0.1cm}

\begin{remark}[Weaker assumptions]\label{WOinInfinStr} In \cite{Gass24B}, we assumed that our metatheory is ZFC and that ZFC is consistent and we used the argument that every standard structure $\S_1^I$ is a model of $\{WO^{n}\}\cup choice_h^{(2)}$ by \cite{Zerm04} to derive the {\rm HPL}-independence of $WO^{n}$. However, it is sufficient that one structure is a model of $\{WO^{n}\}\cup choice_h^{(2)}$. Thus, we can say that the assumptions used here to obtain the result considered once more in Remark \ref{WOinFinStr} are weaker. 
\end{remark}

\vspace{0.1cm}

\section{Why is the basic Fraenkel model a model of HAC?} \label{Section3}

\subsection{A list of ideas for proving HAC}\label{Sect_Idea1}

By Lemma \ref{endlSigma}, the following ideas can be used to prove that any basic Henkin structure $\S$ given by $\S=\S^I_{\sf 0} =_{\rm df}\S(I,\G_{\sf 1} ^I,\I_{\sf 0}^I)$ is a model of $choice_h^{1,1}(H)$. These ideas were also used to prove Proposition 4.14 in \cite{Gass24B}. More precisely, we can say that a generalization of the following ideas also leads to the proof of Proposition 4.14 in \cite{Gass24B}.

\begin{idea}[For proving $choice_h^{1,1}(H)$ considering all $({\cal A}^H_{\xi})_{ \xi \in I}$] \hfill Assuming $\S\models_f\forall x \exists D H( x ,D)$, we will consider a finite choice set $\widetilde{\alpha_0}\subseteq I$ for a partition of $I$ (containing an element of any set of the partition) and take a choice function $\phi_ {|\widetilde{\alpha_0}|-1}$ on the family $({\cal A}^H_{\xi})_{ \xi \in \widetilde{\alpha_0}}$ in order to create a suitable choice function $\varphi$ on the family $({\cal A}^H_{\xi})_{ \xi \in I}$ by extending $\phi_ {|\widetilde{\alpha_0}|-1}$. In doing so, we want to ensure that suitable predicates $\delta_\xi\in {\cal A}_\xi^H$ can be provided in a uniform way by $\delta_{\xi}=\varphi(\xi)$ and that $\sigma_{\varphi}$  (defined in \ref{Conseq} of Lemma \ref{Zhg_CG_sigma}) will be definable in $\S_0$. Therefore, we consider the following points.

\begin{itemize}
\item Let $P$ be a finite stabilizer $\{\nu_1, \ldots,\nu_{|P|}\}$ such that the equation {\rm (\ref{Stabil})} considered in Lemma \ref{EndlichesPfuerH} holds for all $\pi \in \G_{\sf 1} ^I(P)$. $P$ is dependent on $H$, but we omit the index $H$.
\item We decompose $I$ into $q+1$ sets $\{\nu_1\}, \ldots,\{\nu_{q}\}$, and $I\setminus P$ ($q=|P|$) and obtain a $P$-adequate partition.
\item We consider a choice set $\widetilde{\alpha_0}=\{\xi_1, \ldots,\xi_{q+1}\}$ for this partition (later, also denoted by $\widetilde{\alpha^{I,P}_{\mu}}$). 
\item By a choice function $\phi_q$ recursively defined, we obtain $q$ unary predicates $\delta_{\xi_1}=\phi_q({\cal A}_{\xi_1}^H), \ldots , \delta_{\xi_{q+1}}=\phi_q({\cal A}_{\xi_{q+1}}^H)$. 
\item We assume that $\xi_{q+1}\in I\setminus P$.
\item For any $\xi \in I\setminus P$, we try to find a suitable $\pi\in \G_{\sf 1} ^I(P)$ to get $\xi=\pi(\xi_{q+1})$. Then, for this $\xi$, we can choose the predicate $\pi(\phi_q({\cal A}_{\xi_{q+1}}^H))\in {\cal A}_{\xi}^H$.
\begin{itemize}
\item Let $\xi_{q+1}=\mu_1$ for some $\mu_1\in I\setminus P$ (e.g. $\mu_1=\min (\bbbn\setminus P)$ if $I=\bbbn$).
\item Starting with $\mu_1 \in I\setminus P$ and $\phi_q({\cal A}_{\mu_1}^H)$, we will use a formula for determining a permutation $\pi_\xi \in \G_{\sf 1}^I(P)$ with $\xi=\pi_\xi(\mu_1)$ and thus $\pi_\xi(\phi_q({\cal A}_{\mu_1}^H))\in {\cal A}_{\xi}^H$ for this $\xi $. 
\item Then, for $ \delta _{\mu_1}= \phi_q({\cal A}_{\mu_1}^H)$, and any individual $\eta$ given by $ \eta=\pi_{\xi}( \eta_0) $, $\eta_0\in \delta _{\mu_1}$ holds if and only if $\eta\in \pi_\xi(\delta _{\mu_1})$ holds.
\item The formula used will be denoted by $swap_{q+1}$. It will be defined later in order to prove that the predicate $\sigma_{\varphi}$ is definable over $\S_0$.
\end{itemize} 
\end{itemize} 
\end{idea}

\vspace{0.3cm}
\noindent{\bf Summary: A binary predicate $\sigma$ for $H$ and $f$ with $\S_0\models_f\forall x \exists D H( x ,D)$}

\nopagebreak 
\noindent\fbox{\parbox{11.8cm}{

\vspace*{0.1cm}

\begin{tabular}{ll}
$P$& is a stabilizer $\{\xi_1, \ldots,\xi_{q}\}$ with {\rm (\ref{Stabil})} in Lemma \ref{EndlichesPfuerH} for all $\pi \in \G_{\sf 1} ^I(P)$.\\
$\xi_{q+1}\!\!$&$=\mu_1 \in I\setminus P$.\\
$E$&$=\{1,\ldots, q+1\}$.\vspace{0.2cm}\\
$\delta_{\xi_e}$&satisfies $\S_0\models_{f\langle{x\atop \xi_e} {D\atop \delta_{\xi_e}}\rangle} H( x ,D)$
\hfill for $e\in E$.\\
$\xi$&$=\pi_\xi(\xi_{q+1})$ \hfill for $\xi \in I\setminus P$.\\
$\sigma$&$=\bigcup_{e\in E}\{(\xi_e,\eta)\mid \eta\in \delta_{\xi_e}\} \cup \{(\xi,\eta)\mid \xi\in I\setminus P\,\,\&\,\,\eta\in \pi_\xi(\delta_{\xi_{q+1}})\}$.\\
\end{tabular}
}}

\nopagebreak 
\noindent\fbox{\parbox{11.8cm}{
The task is to show that $\sigma$ and $\sigma_\varphi$, respectively, belongs to $\S_0$.
}}

\subsection{$P$-adequate decompositions}\label{Decomp}

Since we want to decompose the domain $I $ --- in particular for $I=\bbbn$ --- by means of unary predicates in a Henkin-Asser structure in ${\sf struc}_{\rm pred}^{({\rm m})}(I)$, we will say that, for predicates $\alpha, \alpha_1, \ldots, \alpha_l \in {\rm pred}_1(I)$ ($l\geq 1$), the set $\{\alpha_1, \ldots, \alpha_l\}$ is a {\em decomposition of $\alpha$} if $\{\widetilde{\alpha_1}, \ldots, \widetilde{\alpha_{l\,}}\}$ is a decomposition of $\widetilde{\alpha}$. A predicate $\alpha$ is {\em empty} if $\widetilde{\alpha}$ is empty. A predicate $\alpha$ is {\em non-empty} if $\widetilde{\alpha}$ is non-empty. A non-empty decomposition that does not contain an empty predicate (or set) is a {\em partition}. We will use many no\-ta\-tions as given in \cite{Gass24B}.

The finite individual supports allow to decompose an individual domain $I$ as follows. Let $P$ be any finite subset of $ I$. If $P$ is not empty, then let $P=\{\nu_1, \ldots, \nu_q\}$ for some individuals $\nu_1, \ldots, \nu_q\in I$ and $q=|P|$, otherwise, let $q=0$. Thus, $|P|=q$. For decomposing the domain $I$, we introduce $q+1$ unary predicates. Let, for all $j\in \{1,\ldots,q\}$,
\[\beta_j=\{\nu_j\} \quad\mbox{ and let }\quad \beta_{q+1}= I\setminus P.\]
The finite set $\{ \beta_1,\ldots, \beta_{q+1} \}$ is a set of mutually exclusive predicates and thus a partition of $I$ which we call the {\em $P$-adequate partition of the individual domain $I$}.

\begin{lemma}[$P$-adequate partition of an individual domain]\label{BetasInHenkM} \hfill For any \linebreak non-empty set $I$ and any finite subset $P$ of $ I$, each predicate $\beta$ in the $P$-adequate partition of $I$ belongs to each Henkin-Asser structure in ${\sf struc}_{\rm pred}^{({\rm m})}(I)$. Moreover, there holds ${\rm sym}_\G(\beta)\supseteq \G(P)$ for each subgroup $\G$ of $\G^I_1$.
\end{lemma}

We follow \cite{Gass24B} and consider the case $n = 1$ and $m=1$. Let $q\geq 0$ and $E= \{1,\ldots,q+1\}$. $E$ depends on $q$ in any case, but we do not write $E_q$, i.e. we omit the index $q$ (which here generally stands for $|P|$). Moreover, let $e$ and $e '$ be any indices in $E$. 

\setcounter{equation}{0}

\begin{lemma}[Formulas describing the $P$-adequate partition of $I$]\label{dieBs} \hfill Let  \linebreak  $\S$ be in ${\sf struc}_{\rm pred}^{({\rm m})}(I)$ for any domain $I$, $f$ be in ${\rm assgn}(\S)$, $P$ be a finite subset of $I$, and $\{ \beta_1,\ldots,\! \beta_{q+1} \}$ be the $P$-adequate partition of $I$. Let $\xi\in I$ and $\pi\in {\G}_{\sf 1}^I(P)$. If $\beta_1,\ldots, \beta_{q+1}$ belong to $\S$, then the statements {\em (\ref{dieBs2})} and {\rm (\ref{dieBs20})} hold for the assignment $\bar f$ given by $\bar f=f\langle{ B_1 \atop\beta_1 } {\cdots\atop \cdots} { B_{q+1}\atop\beta_{q+1}}\rangle$. 
Let $ e \in E$. If $\beta_e$ belongs to $\S$, then {\rm (\ref{dieBs2b})} and {\rm (\ref{dieBs0})} hold for $\bar f=f\langle{ B_e \atop\beta_e} \rangle$. 
\begin{equation}\label{dieBs2}\S\models_{\bar f} \bigvee_{e \in E} (B_ex)\end{equation}
\begin{equation}\label{dieBs20}\S\models_{\bar f}\bigwedge_{e, e'\in E \atop e \not= e ' } (B_ex \to \neg B_{e '}x)\end{equation}
\begin{equation}\label{dieBs2b}\S_{\bar f\langle{x \atop \xi} \rangle}
(B_ex)=\S_{\bar f\langle{ x \atop\pi(\xi)} \rangle} (B_ex)\end{equation}
\begin{equation}\label{dieBs0}\beta_e =\alpha_{\S,B_ex,x,\bar f}\end{equation} 
\end{lemma}

By (\ref{dieBs2b}) and (\ref{dieBs0}) in Lemma \ref{dieBs}, we have ${\rm sym}_\G(\beta_e )\supseteq \G(P)$ for each subgroup $\G$ of $\G^I_1$. The statements (\ref{dieBs2}) and (\ref{dieBs20}) in Lemma \ref{dieBs} imply that the set $\{\beta_e \mid e \in E\}$ of non-empty predicates is a finite decomposition of $I $. 

For the proof of $\S_0\models choice_h^{1,1}(H)$ for $H$ in ${\cal L}^{(2)}_{ x ,D}$, it is useful to take the following questions into account. Let $\xi $ and $\eta$ be individuals in $I$. What happens when we change the values $\xi$ and $\eta$ by a permutation $\pi\in \G_{\sf 1}^I(P)$? What happens when we change some of the values in $I\setminus P$? We know that $\xi=\eta$ holds iff $\pi(\xi)=\pi(\eta)$ holds. Let $\{\beta_1,\ldots, \beta_{q+1}\}$ be the $P$-adequate partition of $I$ and $\beta_{q+1}=I\setminus P$. Moreover, let $\mu_1 \in I \setminus P$ and $P_\mu=\{\mu_1\}$ ($\mu$ is a unary tuple with the only component $\mu_1$).

\begin{lemma}[A choice set for the $P$-adequate partition of $I$]\label{L_P_id}\hfill For the \linebreak predicates $\beta_1,\ldots, \beta_{q+1}$ defined above,
\[\widetilde{\alpha^{I,P}_{\mu }}=_{\rm df} \{\xi_1,\ldots, \xi_{q+1}\}=\{\nu_1,\ldots, \nu_q,\mu_1\}\] is a choice set for the partition $\{\beta_1,\ldots, \beta_{q+1}\}$ of $I$. For all subgroups $\G\subseteq \G_{\sf 1}^I$, ${\rm sym}_{\G}(\{ \xi_e\})\supseteq \G(P\cup P_{\mu })$ holds for all $e\in E$ and thus we have ${\rm sym}_{\G}(\alpha_{\mu }^{I,P})\supseteq \G(P\cup P_{\mu })$. 
\end{lemma}

\subsection{The description of the consequences of swapping two individuals}\label{DefinPerm}

Let $P\in \I_0^I$ and $\mu_1\in I\setminus P$. For $\xi\in I\setminus P$, we want to define a list $\zeta_1,\ldots,\zeta_{s_{\xi}}$ of $s_\xi$ pairwise different individuals in $ I\setminus P$ satisfying $\{\zeta_1,\ldots,\zeta_{s_{\xi}}\}=\{\mu_1,\xi\}$ with $s_\xi\leq 2$ (in case of $n=m=1$) --- which means $ \G^I_{\sf 1}(I\setminus \{\zeta_1,\ldots,\zeta_{s_{\xi}}\})\subset \G^I_{\sf 1}(P)$ --- such that it is possible to define a $\pi_{\xi}\in \G^I_{\sf 1}(I\setminus \{\zeta_1,\ldots,\zeta_{s_{\xi}}\})$ which satisfies $\xi =\pi_{\xi} ( \mu_1)$ and which is describable by a second-order formula. For $\xi\in P$, let $\pi_{\xi}$ be in $\G^I_{\sf 1}(I)$ and $s_{\xi}=0$. In the latter case, $\pi_{\xi}$ is the only mapping in $\G^I_{\sf 1}(I)$, it is an identity bijection.

It will be important for us that the permutation $\pi_{\xi}$ defined as follows is in $\G^I_{\sf 1}(P)$ and that its application to any individual $\eta_0$ that results in $\eta= \pi_{\xi}( \eta_0)$ --- which means $\eta\in \pi_\xi(\delta _{\mu_1})$ for $\eta_0\in \delta _{\mu_1}$ and $ \delta _{\mu_1}= \phi_q({\cal A}_{\mu_1}^H)$ --- can be uniquely described by a finite formula in ${\cal L}^{(2)}_{x_0, y_0, x, y}$ under the assignment $\langle { x_0\atop \mu_1}{ y_0\atop \eta_0}{x \atop \xi}{ y \atop \eta }\rangle $.
If $\xi\in I\setminus P$, let 
\[\zeta_j= \left\{\begin{array}{ll} \mu_1 &\mbox{ if }1= j, \\
\xi &\mbox{ if }1<j \leq s_{\xi} \mbox{ (which means } s_{\xi}=2 \mbox{ and } \xi\not=\mu_1) . 
\end{array}\right.\]
Let $\pi_{\xi}\in \G^I_{\sf 1}(P)$ be a permutation defined step by step by 
\[\pi_{\xi}(\zeta)=
\left\{\begin{array}{ll} \zeta&\mbox{ if } \zeta\in I \mbox{ and } \zeta\not = \zeta_j \mbox{ for all } j\leq s_{\xi},\\
\xi&\mbox{ if } \zeta=\zeta_1 \mbox{ (which means } \zeta=\mu_1),\\
\mu_1 &\mbox{ if } \zeta = \zeta_j \mbox{ for } 1 < j \leq s_{\xi} \mbox{ (which means } s_{\xi}=2 \mbox{ and } \zeta\not=\mu_1) . \end{array}\right.\]
Consequently, we have $\pi_{\mu_1}\in \G_0^I$ and $\pi_\xi(\mu_1)=\xi$ and $\pi_\xi(\xi)=\mu_1$. For the predicates $\gamma$, $\gamma_1$, $\gamma_2$, and $\gamma_3$ given by $\widetilde{\gamma}\subseteq I\setminus \{\mu_1,\xi\}$, $\widetilde{\gamma_1}= \gamma\cup\{\mu_1\}$, $\widetilde{\gamma_2}= \gamma\cup\{\xi\}$, and $\widetilde{\gamma_3}= \gamma\cup\{\mu_1,\xi\}$, we have $\gamma^{\pi_\xi}=\gamma$, $\gamma_1^{\pi_\xi}=\gamma_2$, $\gamma_2^{\pi_\xi}=\gamma_1$, and $\gamma_3^{\pi_\xi}=\gamma_3$.

Let us consider such a permutation for changing (or transforming) pairs $(\mu_1, \eta_0)$ of individuals that belong to $\{\mu_1\}\times \delta_\xi$. The definition of $\pi_\xi\in \G^I_{\sf 1}(P)$ implies that the relationship between the pairs $( \mu_1, \eta_0)$ and $( \xi ,\pi_{\xi}( \eta_0))$ can be described by one formula $swap_{q+1}$ --- for the fixed $n,m= 1$ --- under an assignment of the form $\langle {x_0\atop \mu_1}{ y_0\atop \eta_0}{x \atop \xi}{ y \atop \eta }\rangle $.

We will use $x_0$, $x$, $y_0$, and in particular $ y $ for describing $ \eta =\pi_{\xi} ( \eta_0)$. 

\vspace{0.2cm}

\noindent If the variable $x_0$ stands for an individual $\xi_e$ in $P$ with $e\in E\setminus \{q+1\}$, then let $swap_e(x_0, y_0, x, y)$ stand for $x = x_0 \land y = y_0 $. If $x_0$ stands for an individual in $I\setminus P$, let $swap_{q+1}(x_0, y_0, x, y)$ stand for
\[(x = x_0 \to y = y_0 ) \,\,\land\,\, ( x \not = x_0 \to\hspace*{6.7cm}\]
\[(( y_0 \not=x_0\land y_0 \not=x) \to y=y_0)\land (y_0 =x_0 \to y=x)\land (y_0 =x \to y = x_0)).\]
Let us remark that, (in particular here) for $n=m=1$, it would be possible to take the latter formula for all $e\in E$. Hence, we could also omit the index $e$. 

\begin{lemma}[A formula for permuting individuals]\label{Perm_Swap} Let $\S$ be a structure in ${\sf struc}_{\rm pred}^{({\rm m})}(I)$ and $f$ be in ${\rm assgn}(\S)$. Let $P$ be a finite subset of $I$, let $q=|P|$, and let $e\in E$. Let $\{\beta_1,\ldots, \beta_{q+1}\}$ be the $P$-adequate partition of $I$, $\beta_{q+1}=I\setminus P$, and $\xi_e$ be in $\widetilde{\beta_e}$. Moreover, let $\xi$ be in $\widetilde{\beta_e}$, $\pi_{\xi}\in \G^I_{\sf 1}(P)$ be defined as above, and $ \eta_0\in I$. Then,
\[\S\models_{f'} swap_e(x_0, y_0, x, y)\] holds if and only if $f'=f\langle { x_0\atop \xi_e}{ y_0\atop \eta_0}{x \atop \xi}{ y \atop\pi_{\xi}( \eta_0) }\rangle $. 
\end{lemma}
The latter means that $f'=f\langle { x_0\atop \mu_1}{ y_0\atop \eta_0}{x \atop \xi}{ y \atop \eta }\rangle $ holds for $ \eta =\pi_{\xi} ( \eta_0)$.

Now, the preparation for giving the details of a proof that the 1-1 Ackermann axioms are true in the basic model is closed.

\vspace*{0.3cm}
\newcounter{liass}
\begin{table}[t]
\noindent\framebox{$\!\!\!\!$\parbox{12.05cm}{\begin{list}{$\!\!\!\!\!\!$}
{\usecounter{liass}}
\item\vspace*{-3pt}$\S$ \,\,is a predicate structure $(J_n)_{n\geq 1}$ in ${\sf struc}_{\rm pred}^{({\rm m})}(I)$ for $I=\bbbn$. 
\item\vspace*{-3pt}$P$ \,\,is a finite subset of $I$.
\item\vspace*{-3pt} $q \,\, =|P|$.
\item\vspace*{-3pt}$E=\{1,\ldots,q+1\} $.
\item\vspace*{-3pt} $e \,\, \in E$.
\item\vspace*{-3pt} ${\rm idx}_{q+1} \,\, =\{1\}$.
\item\vspace*{-3pt} ${\rm idx}_{\nu} \,\, =\emptyset$ for $\nu\in\{1,\ldots,q\}$.
\item\vspace*{-3pt} ${\sf K}_{{\rm idx}_ {e }}$ \,\,contains all partitions of ${\rm idx}_ {e }$.
\item\vspace*{-3pt} ${\cal K} \,\, \in {\sf K}_{ {\rm idx}_ {e }}$. 
\item\vspace*{-3pt} $\vec {\cal K} \,\,$ has one component $\{1\}$ if $e=q+1$ and thus ${\cal K}=\{\{1\}\}$.
\item\vspace*{-3pt} $F_{\cal K}(x) \,\, = x =x$.
\item\vspace*{-3pt} $\alpha_{I,{\cal K}}(\xi) \,\, = \S_{\langle{x \atop \xi} \rangle} (F_{\cal K}(x))$ (and thus $\widetilde{\alpha_{I,{\cal K}}}=I$).
\item\vspace*{-3pt} $\{\beta_1,\ldots, \beta_{q+1}\}$\,\,is the $P$-adequate partition of $I$ with $\beta_{q+1}=I\setminus P$.
\item\vspace*{-3pt} $F_{e, {\cal K}}(x) \,\, = B_{e}x \land F_{\cal K}(x)$ can be replaced by $ B_{e}x$. 
\item\vspace*{-3pt} $\alpha_{e, {\cal K}}^{I,P}(\xi) \,\, = \S_{\langle{x \atop \xi} { B_1 \atop\beta_1} {\cdots\atop\cdots} { B_{q+1}\atop\beta_{q+1}}\rangle} (F_{e , {\cal K}}(x))$ if $\beta_1,\ldots,\beta_{q+1}\in J_1$.
\item\vspace*{-3pt} $\mu$\,\, is a tuple $(\mu_1)\in \beta_{q+1}$.
\item\vspace*{-3pt} $P_{\mu} \,\, =\{\mu_1\}$.
\item\vspace*{-3pt} $\xi_{\mu,e , {\cal K}}^{I,P}$\,\, (shortly, $\xi_{e , {\cal K}}$) is here $\xi_e$ with $\xi_e\in \{\nu_1,\ldots, \nu_q,\mu_1\}$.
\item\label{ASS11}\vspace*{-3pt}
$\alpha^{I,P}_{\mu}(\xi) \,\, =true$ iff $\xi\in \{\xi_{\mu,e , {\cal K}}^{I,P} \mid e\in E\}$.
\end{list}}}
\caption{Partitions of $\bbbn$ and some special notations for $n=1$ and $m=1$.}\label{overview_nota_1}
\end{table}

\subsection{The 1-1 Ackermann axioms hold in the basic model}\label{HACinBasic}

Let us recall that $\S_0$ is the Henkin-Asser structure $\S(\bbbn,\G_{\sf 1} ^\bbbn,\I_{\sf 0}^\bbbn)$ which we also call the basic Freankel model of second order. Now, we will confirm
\[\S_0\models choice_h^{n,m}(H)\]
 for $n=1$ and $m=1$ and $H$ in ${\cal L}_{ x ,D}^{(2)}$.

\vspace{0.3cm}

\noindent {\bf A proof for $\S_0\models choice_h^{1,1}$.} 

Let $\S$ be the basic structure $\S(\bbbn,\G_{\sf 1} ^\bbbn,\I_{\sf 0}^\bbbn)$ with the domains $J_0=\bbbn$ and $J_1=J_1(\bbbn,\G_{\sf 1} ^\bbbn,\I_{\sf 0}^\bbbn)$ and let $\G=\G_{\sf 1} ^\bbbn$. Let $f\in {\rm assgn}(\S)$ and $H( x ,D)$ be a formula in ${\cal L}_{ x ,D}^{(2)}$ such that $\S\models_f \forall x \exists D H( x ,D) $ holds. Let $P$ be a finite subset of $\bbbn$ such that the equation {\rm (\ref{Stabil})} in Lemma \ref{EndlichesPfuerH} holds for all $\xi \in \bbbn$, all predicates $\delta\in J_{1}$, and all $\pi\in \G(P)$. Such a stabilizer $P$ exists by Lemma \ref{EndlichesPfuerH}. Let $\{ \beta_1,\ldots,\beta_{q+1} \}$ be the $P$-adequate partition of $I$ and $\beta_{q+1}=\bbbn \setminus P$. It follows from Lemma \ref{BetasInHenkM} that the predicates $\beta_1,\ldots,\beta_{q+1}$ belong to $J_1$. In the following, let $\mu=\mu_1=\min (\bbbn \setminus P)$ and let $P_{\mu }=\{\mu_1\}$. Let $B_1, \ldots, B_{q+1}$ be new variables for $1$-ary predicates that do not occur in $H( x ,D)$. Then, according to Lemma \ref{dieBs}, for $\bar f=f\langle { B_1 \atop\beta_1 } {\cdots\atop\cdots} { B_{q+1}\atop\beta_{q+1}}\rangle $, $\xi \in \bbbn $, $\delta \in J_1$, and each $\pi\in\G(P)$, we have 
\[\!\!\!\! \!\!\S_{\bar f\langle{x \atop \xi}{D\atop\delta}\rangle} (B_e x \land H( x ,D))=\S_{\bar f\langle{ x \atop \pi(\xi)}{D\atop\delta^\pi}\rangle} (B_e x\land H( x ,D))\] for each $e\in E$. Let $\alpha_0$ be the unary predicate $\alpha^{\bbbn,P}_{\mu }$ defined as $\alpha^{I,P}_{\mu }$ for any individual domain $I$ and considered in Lemma \ref{L_P_id}. Thus, $\widetilde{\alpha_0}$ is the choice set $\{\nu_1,\ldots, \nu_q, \mu_1\}$ considered in Lemma \ref{L_P_id}. $\widetilde{\alpha_0}$ is finite and, thus, $\alpha_0$ belongs to $J_1$. By assumption, $\S_{f\langle{x \atop \xi}\rangle} (\exists D H( x ,D))$ is $true$ for all $\xi \in \bbbn $. Hence, for any $\xi \in \alpha_0$ there is at least one predicate $\delta \in J_1$ such that 
\[\S\models_{f\langle{x \atop \xi} {D\atop\delta }\rangle} H( x ,D).\]
Consequently, we have a finite family $({\cal A}^H_{\xi})_{ \xi \in \alpha_0} $ of non-empty sets $ {\cal A}^H_{\xi}$ with $ {\cal A}^H_{\xi}=\{\delta\in J_1\mid \S\models_{f\langle{ x \atop \xi} {D\atop\delta }\rangle}H( x ,D)\}$ and thus
\[\S\models_{f\langle {A\atop\alpha_0}\rangle} \forall x \exists D(A x \to H( x ,D)).\]
 This implies, by using a proof of Lemma \ref{endlSigma}, the existence of a choice function $\phi_s$ that can be recursively defined and then used, for $s=q$, in the definition of the $(1+1)$-ary predicate $\sigma_0\in J_{1+1}$ satisfying
 \[\widetilde{\sigma_0} =\{(\xi , \eta)\in \bbbn^{1+1}\mid \xi \in \alpha_0 \,\,\&\,\, \eta \in\phi_s ( {\cal A}^H_{\xi})\}.\]
For each $\xi \in \alpha_0$, let the predicate $\delta_{\xi}$ be defined by $\delta_{\xi}=\phi_s( {\cal A}^H_{\xi})$. 
Then, for all $ e \in E $, $\delta_{ \xi_e} $ is --- by the definition of $ {\cal A}^H_{ \xi_e}$ --- in $\S$ and, therefore, there is some finite $P_{\delta_{ \xi_e}}\subseteq\bbbn$ such that ${\rm sym}_\G(\delta_{ \xi_e})\supseteq\G(P_{\delta_{ \xi_e}})$. Let 
\[\sigma_e^{(0)}=\{( \xi_e, \eta )\in \bbbn^{1+1}\mid \eta \in \delta_{ \xi_e}\}.\] 
Since ${\rm sym}_{\G}(\sigma_e^{(0)})\supseteq \G(P\cup P_{\mu }\cup P_{\delta_{ \xi_e}})$ follows from ${\rm sym}_{\G}(\{ \xi_e\})\supseteq \G(P\cup P_{\mu })$ (cf.\,\,Lemma \ref{L_P_id}), the predicate $\sigma_e^{(0)}$ belongs to $J_{1+1}$. Consequently, the predicate $\sigma_0$ can be described by \[\sigma_0= \bigcup_{e \in E}\, \sigma_e^{(0)}.\]
$\sigma_0$ is thus the predicate $\alpha _{\S,G_0 (x , y), (x, y), \tilde f}$ given by $ G_0 (x , y)= \bigvee_{i=1,\ldots,q}R_i x y $ and $\tilde f= f\langle{R_1\atop \rho_1}{\cdots\atop\cdots }{R_{q}\atop \rho_{q}}\rangle$ if we assume that the set $ \{\rho_1,\ldots,\rho_{q}\}$ is equal to $ \{\sigma_e^{(0)}\mid e \in E \} $. This means that $\sigma_0$ belongs to $J_{1+1}$. Moreover, the relationship ${\rm sym}_{\G}(\sigma_0)\supseteq \bigcap_{e \in E}\, \G(P\cup P_{\mu }\cup P_{\delta_{ \xi_e}})$ is satisfied. Let $P_0$ be the finite set $ \bigcup_{e \in E}\, P_{\delta_{ \xi_e}}$. Thus, we have ${\rm sym}_{\G}(\sigma_0)\supseteq\G(P\cup P_{\mu }\cup P_0 )$ and the set $P\cup P_{\mu }\cup P_0$ is therefore a finite individual support for $\sigma_0$. By definition of $\sigma_0$, there holds
 \[\S\models_{f\langle{A\atop\alpha_0}{S\atop \sigma_0} \rangle}\forall x \exists D ( A x \to \forall y (D y \leftrightarrow S x y ) \land H( x , D )).\] 
Now, let $ e \in E $. Let $\pi_{\xi}$ be the permutation defined for all $\xi \in I$ in Section \ref{DefinPerm} and let 
\[\sigma_e = \{(\xi , \eta)\in \bbbn^{1+1}\mid \eta \in (\delta_{ \xi_e})^{\pi_{\xi}}\}\] 
which means that the tuple $( \xi , \eta )$ belongs to $\sigma_e $ iff there is an $ \eta_0\in \delta_{ \xi_e }$ such that $ \eta =\pi_{\xi}( \eta_0)$. Let $\bar f=f\langle{ B_1 \atop\beta_1 } {\cdots\atop\cdots} { B_{q+1}\atop\beta_{q+1}} \rangle $ and $f'=\bar f\langle {D\atop \delta_{ \xi_e}} { x_0\atop \xi_e} \rangle $. For 
\[G_e(x, y)= B_e( x )\land \exists y_0(D y_0\land swap_e(x_0, y_0, x, y)),\]
we get $\sigma_e =\alpha _{\S,G_e (x, y), (x, y), f'}$ according to Lemma \ref{Perm_Swap}. Consequently, $\sigma_e \in J_{1+1}$ and ${\rm sym}_{\G}(\sigma_e )\supseteq \G(P\cup P_{\mu }\cup P_{\delta_{\xi_e} })$ hold. Let $f''= f'\langle{S_0\atop \sigma_0} \rangle$ and $\bar f'= \bar f\langle{S_0\atop \sigma_0} \rangle$. Because of $\S\models_{f''} D y \leftrightarrow S_0x_0y$, we have the equation $\S_{f'} ( G_e(x, y))= \S_{\bar f'}( G_e '(x, y) )$ for
\[G_e '(x, y)=B_e( x ) \land\exists x_0 \exists y_0 (B_e(x_0) \land S_0x_0y_0 \land swap_e(x_0, y_0, x, y))\] 
which implies
\[\sigma_e =\alpha _{\S,G_e(x, y), (x, y), f'} =\alpha _{\S,G'_e (x, y), (x, y), \bar f'}.\]
Let $\xi$ satisfy $\S\models_{\bar f\langle{x \atop \xi}\rangle} B_e( x )$ and let $\delta_{\xi}$ be the predicate $(\delta_{ \xi_e})^{\pi_{\xi}} $. 
Then, there holds $\{ \eta \mid (\xi , \eta) \in\sigma_e \}= \delta_{\xi}$. Since ${\rm sym}_{\G}(\delta_{ \xi_e})\supseteq \G(P_{\delta_{ \xi_e}})$ implies ${\rm sym}_{\G}(\delta_{\xi})\supseteq \G(\pi_{\xi}(P_{\delta_{ \xi_e}}))$, $\delta_{\xi}$ belongs to $\S$, too (where $\pi_{\xi}(P_{\delta_{ \xi_e}})=\{\pi_{\xi}(\eta_0)\mid \eta_0\in P_{\delta_{ \xi_e}}\}$). By Lemma \ref{EndlichesPfuerH} and by the definition of the stabilizer $P$, for any $\pi\in \G(P)$, we have
\[\S_{f\langle { x \atop \xi_e} {D\atop\delta_{ \xi_e}}\rangle}(H( x ,D))=\S_{f\langle { x \atop \pi( \xi_e )} {D\atop(\delta_{ \xi_e})^{\pi}}\rangle}( H( x ,D))=\S_{f\langle{x \atop \xi} {D\atop\delta_{\xi}}\rangle}(H( x ,D))\] 
which means $\delta_{\xi} $ is in ${\cal A}_{\xi}^H$ (defined in Lemma \ref {Zhg_CG_sigma}) and
\[\S\models_ {\bar f\langle {S\atop \sigma_e }\rangle}\forall x (B_e( x )\to \exists D (\forall y (D y \leftrightarrow S x y ) \land H( x , D ))).\] 
On the other hand, the predicate
 \[\sigma= \bigcup_{e \in E}\sigma_e\] 
that can be also described by $\sigma =\alpha _{\S,G(x , y), (x, y), \bar f\langle {S_0\atop \sigma_0 }\rangle}$ and 
\[G(x , y)= \bigvee_{ e \in E}G_e '(x , y)\]
 --- which implies ${\rm sym}_{\G}(\sigma)\supseteq \G(P\cup P_{\mu }\cup P_0)$ --- belongs to $J_{1+1}$ and satisfies
\[\S\models_{f\langle{S\atop \sigma} \rangle}\forall x \exists D ( \forall y (Dy \leftrightarrow S xy) \land H( x , D )). \] 
This means that the consequent of $choice_h^{1,1}(H)$ follows from the antecedent of $choice_h^{1,1}(H)$. 
\hfill \qed

\vspace{0.2cm}

The proof confirms that we have $\S_0 \models {\rm HAC}$ and that there is a Henkin-Asser structure that is a model of 
\[^{h}ax^{(2)} \cup choice_h^{1,1} \cup \{\neg WO^1 \}.\] 

\noindent The latter is a special case of a result proved in \cite{Gass24B}.

\section{HAC, TR, and a finite union of Fraenkel models}\label{Section4}

We can prove that the Henkin-Asser structures $k\S_0$ are, for $k\geq 2$, models of $choice_h^{(2)}\cup \{\neg TR^n\}$ ($n\geq 1$). But we will only prove that a Henkin-Asser structure of this form is a model of ${\rm HAC}\cup \{\neg TR^1\}$. A generalization of the following proof by a combination of the ideas used here with the ideas, that are used in the proof of $\S_0\models choice_h^{n,m}$ (see \cite{Gass24B}) in order to consider the effects of permutations on $n$- or $m$-tuples (some of whose components can be equal), leads to a proof of $k\S_0\models choice_h^{n,m}$. The combination of the ideas is possible since all unary predicates in $k\S_0$ result from a union $\widetilde{\alpha_1} \cup\cdots \cup \widetilde{\alpha_k}$ given by $k$ finite or co-finite predicates $\alpha_1\subseteq \{1\}\times\bbbn ,\ldots, \alpha_k\subseteq \{k\}\times\bbbn $ belonging to $k$ pairwise disjoint copies of the basic Fraenkel model II. Let us remark that the notion {\sf co-finite} is here related to a basic set of the form $\{j\}\times\bbbn$ (with $j\in \{1,\ldots,k\}$). This implies that the proof of $choice_h^{n,m}$ ($n,m\geq 1$) can be done by using finite refinements of $P$-and-${\rm id}$ adequate partitions of $I^n$ ($n\geq 1$).

\subsection{The characterization of a finite union of Henkin structures}\label{Chara}

Let us recall some definitions. Let $I$ be a domain of individuals, $\G$ be a subgroup of $\G_{\sf 1} ^I$, and $\G^*= \{ \pi^* \mid \pi\in \G\}$ where, for each $n\geq 1$, $(\pi^*(\alpha))(\mbm{\xi})=\alpha^\pi(\mbm{\xi})$ holds for all $\alpha\in {\rm pred}_n(I)$ and all $\mbm{\xi}\in I^n$. Consequently, $\G^*$ is a subgroup of $(\G_{\sf 1} ^I)^*$. $\G^*$ and $(\G_{\sf 1} ^I)^*$ are groups that act on $I\cup {\rm pred}(I)$ from the left. Let $\S$ be any predicate structure $(J_{n})_{n\geq 0}$ in ${\sf struc}_{\rm pred}^{({\rm m})}(I)$ and let $\G^*|_{\S}$ be formally given by
\[\begin{array}{ll}\G^*|_{\S}=&\{ g:\bigcup _{n\geq 0}J_n\to I\cup {\rm pred}(I) \mid\\&(\exists \pi \in \G)((\forall \xi \in J_0) (g(\xi)=\pi(\xi))\,\, \&\,\,\,(\forall \alpha \in\bigcup _{n\geq 1}J_n) (g(\alpha)=\alpha^{\pi}))\}.\end{array}\]
It is the set of all functions that are the permutations in $\G^*$ restricted to the subset $\bigcup_{n\geq 0}J_n\subseteq I\cup {\rm pred}(I)$. For any $\pi\in \G ^I_{\sf 1}$ and any subset $J\subseteq \bigcup _{n\geq 0}J_n$, let $\pi^*|_{J}$ be the restriction of $\pi^*\in (\G ^I_{\sf 1})^*$ to $J$ and let $ \pi^ *|_\S=\pi ^*|_{\bigcup _{n\geq 0}J_n}$. By definition, we thus have $ \pi^ *|_\S=\bigcup _{n\geq 0} \pi^ *|_{J_n}$. If the extension $ \pi^ *|_\S\in \G^*|_{\S}$ is in $ {\rm perm}(\bigcup _{n\geq 0}J_n)$, then $( \pi^ *|_{J_n})_{n\geq 0}$ is in ${\rm auto}(\S)$. Let $ {\rm auto}_\G^*(\S)$ and $ {\rm auto}^*(\S)$ be defined by $ {\rm auto}_\G^*(\S)= \G^*|_{\S}\cap {\rm perm}(\bigcup_{n\geq 0}J_{n})$ and $ {\rm auto}^*(\S)= {\rm auto}_{\G ^I_{\sf 1}} ^*(\S)$. 

For $k\S_0=\S(I,\G,\I_0^I)$ given above, $\G$ has the property
\setcounter{equation}{0}
\begin{equation}\label{gruppe}\G=\{\pi\in \G_1^I\mid (\forall j\in \{1,\ldots,k\} )(\pi^ *|_{\S_0^{\{j\}\times \bbbn}}\in {\rm auto}^*(\S_0^{\{j\}\times \bbbn}) )\}.\end{equation} 

\subsection{The current list of ideas for proving HAC}\label{Sect_Idea2}

First, we prove that $k\S_0\models choice_h^{1,1}(H)$ holds for $H\in{\cal L}^{(2)}_{x,D}$. Since $1\S_0$ is the basic Henkin structure $\S_0^{\{1\}\times \bbbn}$ and thus isomorphic to the basic Fraenkel model II, let us consider only the case $k\geq 2$. Now, we again have $I=\{1,\ldots,k\}\times \bbbn$, $I_j=\{j\}\times \bbbn$ for all $j\in\{1,\ldots,k\}$, and the group $\G$ given by (\ref{gruppe}) in Section \ref{Chara}.

\vspace*{0.3cm}

\begin{table}[h]
\noindent\framebox{$\!\!\!\!$\parbox{12.05cm}{\begin{list}{$\!\!\!\!\!\!$}
{\usecounter{liass}}
\item\vspace*{-3pt}$\S$ \,\,is a predicate structure $(J_n)_{n\geq 1}$ in ${\sf struc}_{\rm pred}^{({\rm m})}(I)$ for $I=\{1,\ldots,k\}\times \bbbn$. 
\item\vspace*{-3pt}$P$ \,\,is a finite subset of $I$.
\item\vspace*{-3pt} $q \,\, =|P|$.
\item\vspace*{-3pt} $\{\beta_1,\ldots, \beta_{q+1}\}$\,\, is the $P$-adequate partition of $I$ with $\beta_{q+1}=I\setminus P$.
\item\vspace*{-3pt} $\beta_{(j,i)}= I_j \cap \beta_i$ \,\, for all $j\in \{1,\ldots,k\}$ and $i\in \{1,\ldots, q+1\}$.
\item\vspace*{-3pt}$\{ \beta_{(1,1)}, \ldots, \beta_{(k,q+1)} \}\setminus \{\emptyset\}$ is a refinement of $\{\beta_1,\ldots, \beta_{q+1}\}$.
\item\vspace*{-3pt} $E = \bigcup_{j\in \{1,\ldots,k\} \atop i\in \{1,\ldots,q+1\}} \{(j,i)\mid \beta_{(j,i)}\not = \emptyset\}$ \,\, (E is dependent on $q=|P|$).
\item\vspace*{-3pt} $e \,\, \in E$.
\item\vspace*{-3pt} $\alpha_{e }^{I,P}(\xi) \,\, = \,\, \S_{\langle{x \atop \xi} { B_{(j,i)} \atop\beta_{(j,i)}}\rangle} (B_{e}(x))$ \,\, for $e=(j,i)$ and $\beta_{(j,i)}\in J_1$.
\item\vspace*{-3pt} $\mbm{\mu}$\,\, is a tuple $(\mu_1,\ldots,\mu_k)$ in $ \beta_{(1,q+1)}\times \cdots\times \beta_{(k,q+1)}$.
\item\vspace*{-3pt} $P_{\mbmss{\mu}} \,\, =\{\mu_1,\ldots,\mu_k\}$.
\item\vspace*{-3pt} $\xi_{\mu,e , {\cal K}}^{I,P}$\,\, (shortly, $\xi_{e , {\cal K}}$) is here $\xi_e$ with $\xi_e\in \{\nu_{(1,1)},\ldots, \nu_{(k,q)},\mu_1,\ldots, \mu_k\}$.
\item\vspace*{-3pt}
$\alpha^{I,P}_{\mbmss{\mu}}(\xi) \,\, =true$ iff $\xi\in \{\xi_{e } \mid e\in E\}$.
\end{list}}}
\caption{Partitions of $\{1,\ldots,k\}\times \bbbn$ and special notations for $n, m=1$.}\label{overview_nota_2}
\end{table}

\begin{idea}[For proving $choice_h^{1,1}(H)$ for any families $({\cal A}^H_{\xi})_{ \xi \in I}$]\label{Idee2} \hfill
Assuming \linebreak that $k\S_0$ is a model of $\forall x \exists D H( x ,D)$, we will consider a finite choice set $\widetilde{\alpha_0}$ for a partition of $I$ that can be refined by using $\widetilde{\alpha_0}=\bigcup_{j \in\{1,\ldots,k\}} \widetilde{\alpha_{j,0}}$ with $\widetilde{\alpha_{j,0}}\subseteq I_j$ in order to get a choice function $\varphi$ on the family $({\cal A}^H_{\xi})_{ \xi \in I}$ (with $\varphi(\xi) \in {\cal A}^H_{\xi}$ and $\delta_\xi=\varphi(\xi)$) by using choice functions $\phi_ {j, |\alpha_{j,0}|-1}$ on the finite families $({\cal A}^H_{\xi})_{ \xi \in \widetilde{\alpha_{j,0}}}$.

\begin{itemize}
\item Let $P$ be a finite stabilizer $\{\nu_1, \ldots,\nu_{|P|}\}$ 
 such that the equation {\rm (\ref{Stabil})} considered in Lemma \ref{EndlichesPfuerH} holds for all $\pi \in \G(P)$.
\item We decompose $I$ into $q+1$ sets given by $\beta_1=\{\nu_1\}, \ldots,\beta_q=\{\nu_{q}\}$, and $\beta_{q+1}=I\setminus P$ ($q=|P|$). 
\item Let $j\in \{1,\ldots,k\}$ and $P_j=P\cap I_j$.
\item We decompose $I_j$ into sets given by $\beta_{(j,i)}= I_j \cap \beta_i$ for all $i\in \{1,\ldots, q+1\}$ to obtain a partition of $I_j$ by $\bigcup_{i\in E_j}\{ \beta_{(j,i)}\}= \{ \beta_{(j,1)}, \ldots, \beta_{(j,q+1)} \}\setminus \{\emptyset\}$.

This means that we have $E_j=\{i\mid \beta_{(j,i)}\not=\emptyset \}$. 
\item We consider a choice set $\widetilde{\alpha_{j,0}}=\{\nu_{(j,i)}\mid i\in E_j\}$ with $\{\nu_{(j,i)}\}= \widetilde{\alpha_{j,0}} \cap \beta_{(j,i)}$ for all $i\in E_j$ (where we can use that $\beta_{(j,i)}$ is non-empty). 

We denote $\bigcup_{j\in \{1,\ldots,k\}} \alpha_{j,0}$ by  $\alpha_0$ (and could also denote it by $\alpha^{I,P}_{\mbmss{\mu}}$). 
\item By using an arbitrary choice function $\phi_{j,s_j-1}$ on $({\cal A}^H_{\xi})_{ \xi \in \widetilde {\alpha_{j,0}}}$ with $s_j=|\widetilde {\alpha_{j,0}|}$, we can define unary predicates $\delta_{\nu_{(j,i)}}=\phi_{j,s_j-1}({\cal A}_{\nu_{(j,i)}}^H)$ for $i\in E_j$. 
\item For any $\xi \in \widetilde{\beta_{j,q+1}}$, we try to find a suitable $\pi_\xi\in \G(P)$ to get $\xi=\pi_\xi(\nu_{(j,q+1)})$. Then, we can take the predicate $\pi_\xi(\phi_{j,s_j-1}({\cal A}_{\nu_{(j,q+1)}}^H))$ that belongs to $ {\cal A}_{\xi}^H$ (and could be denoted by $\delta_\xi$).
\begin{itemize}
\item Let $\nu_{(j,q+1)}=\mu_j=(j,n_j)$ with $n_j= \min \{n\in \bbbn\mid (j,n)\in I_j\setminus P\}$.
\item Starting with $\mu_j\in I_j\setminus P$ and $\phi_{j,s_j-1}({\cal A}_{\mu_j}^H)$, we will use a formula for determining a $\pi_\xi$ with $\xi=\pi_\xi(\mu_j)$ and thus $\pi_\xi(\phi_{j,s_j-1}({\cal A}_{\mu_j}^H))\in {\cal A}_{\xi}^H$ for this $\xi $. 
\item Then, for $ \delta _{\mu_j}= \phi_{j,s_j-1}({\cal A}_{\mu_j}^H)$ and any individual $\eta\in I$ given by $ \eta=\pi_{\xi}( \eta_0) $, $\eta_0\in \delta _{\mu_j}$ implies $\eta\in \pi_\xi(\delta _{\mu_j})$.
\item This formula $swap_{(j,q+1)}$ can be defined as above (if we assume that its application will be restricted to $\xi\in \beta_{(j,q+1)}$). 
\end{itemize} 
\end{itemize} 
\end{idea}

\vspace{0.3cm}
\noindent{\bf Summary: A binary predicate $\sigma$ for $H$ and $f$ with $k\S_0\models_f\forall x \exists D H( x ,D)$}

\nopagebreak 
\noindent\fbox{\parbox{11.8cm}{

\vspace*{0.1cm}

\begin{tabular}{ll}
$P$ & is a stabilizer $\{\xi_1, \ldots,\xi_{q}\}$ with {\rm (\ref{Stabil})} in Lemma \ref{EndlichesPfuerH} for all $\pi\! \in\! \G(P)$.\\
$\beta_{e}$&$=I_j\cap \{\xi_i\}$ \hfill for $e=(j,i)\in \{1,\ldots,k\}\times \{1,\ldots,q\}$. \\
$\xi_e$&$=\xi_i$ \hfill for $e=(j,i)\in \{1,\ldots,k\}\times \{1,\ldots,q\}$ if $\beta_{e}\not=\emptyset $. \\
$\beta_{(j,q+1)}\!\!\!\!\!$&$=I_j\setminus P$ \hfill for $j \in \{1,\ldots,k\}$. \\
$\xi_{(j,q+1)}\!\!\!\!\!$&$=\mu_j\in I_j\setminus P$ \hfill for $j\in \{1,\ldots,k\}$.\\
$E$&$=\{(j,i)\in \{1,\ldots,k\} \times \{1,\ldots,q+1\}\mid \beta_{(j,i)}\not= \emptyset\}$.\vspace{0.2cm}\\
$\delta_{\xi_e}$& satisfies $k\S_0\models_{f\langle{x\atop \xi_e} {D\atop \delta_{\xi_e}}\rangle} H( x ,D)$
\hfill for $e\in E$.\\
$\xi$&$=\pi_\xi(\xi_{(j,q+1)})$ \hfill for $j\in \{1,\ldots,k\}$ and $\xi \in \beta_{(j,q+1)}$.\\
$\sigma$&$=\bigcup_{e\in E}\{(\xi_e,\eta)\mid \eta\in \delta_{\xi_e}\} $\\ &\quad$\cup \bigcup_{j\in \{1,\ldots,k\}} \{(\xi,\eta)\mid \xi\in \beta_{(j,q+1)}\,\,\&\,\,\eta\in \pi_\xi(\delta_{\xi_{(j,q+1)}})\}$.\\
\end{tabular}
}}

\noindent\fbox{\parbox{11.8cm}{
The task is to show that $\sigma$ belongs to $k\S_0$.
}}

\subsection{The 1-1 Ackermann axioms hold in a union of Fraenkel models}

\begin{proposition}\label{ChoiceInS2} For any $k\geq 1$ and any $H$ in ${\cal L}_{ x ,D}^{(2)}$, there holds
\[k\S_0\models choice_h^{1,1}(H).\]
\end{proposition}

\noindent {\bf Proof.} Let $k\geq 2$ and $\S$ be the structure $k\S_0$ with the domains $J_0$ and $J_1$ given by $J_0=I=_{\rm df}\{1,\ldots,k\}\times \bbbn$ and $J_1=J_1(I,\G,\I_{\sf 0}^I)$ for the group $\G$ given by (\ref{gruppe}) in Section \ref{Chara}. Let $f\in {\rm assgn}(\S)$ and $H( x ,D)$ be a formula in ${\cal L}_{ x ,D}^{(2)}$ such that $\S\models_f \forall x \exists D H( x ,D) $ holds. Let $P$ be a finite stabilizer satisfying the equation {\rm (\ref{Stabil})} in Lemma \ref{EndlichesPfuerH} for all $\xi \in I$, all predicates $\delta\in J_{1}$, and all $\pi\in \G(P)$. Let $\{ \beta_1,\ldots,\beta_{q+1} \}$ be the $P$-adequate partition of $I$ with $\beta_{q+1}=I \setminus P$ and let $\beta_{(j,i)}= I_j \cap \beta_i$ for all $j\in \{1,\ldots,k\}$ and $i\in \{1,\ldots, q+1\}$. Please note that the set $\{ \beta_{(1,1)}, \ldots, \beta_{(k,q+1)} \}\setminus \{\emptyset\}$ is a refinement of the $P$-adequate partition of $I$. The fact that this refinement is a finite partition of $I$ implies the existence of a finite choice set $\alpha_0$ for this partition and the existence of a useful choice function $\phi_s$ ($s=|\alpha_0|-1$) definable analogously as in a proof of Lemma \ref{endlSigma}. Let us give more details. Let $ E = \bigcup_{j\in \{1,\ldots,k\} \atop i\in \{1,\ldots,q+1\}} \{(j,i)\mid \beta_{(j,i)}\not = \emptyset\}$ and $E_{q+1}=\{1,\ldots,k\} \times \{q+1\}$. For all $e\in E$, $\beta_e$ belongs to $J_1$ because of ${\rm sym}_{\G}(\beta_e)\supseteq \G(P)$. According to the definitions, $\beta_{(j,q+1)}$ is not empty for all $j\in \{1,\ldots,k\}$. Therefore, we can choose and we choose the individual $\mu_j$ in $\beta_{(j,q+1)}$ defined by $\mu_j=(j,n_j)$ with $n_j= \min \{n\in \bbbn\mid (j,n)\in I_j\setminus P\}$. For any $e\in E\setminus E_{q+1}$, let $\nu_e$ be the only individual in $\beta_e$ which means $\beta_e=\{\nu_e\}$. Let $\alpha_0$ be the predicate given by $\widetilde{\alpha_0}=\{\nu_e\mid e\in E\setminus E_{q+1}\}\cup \{\mu_1,\ldots,\mu_k\} $ and let $P_{\mbmss{ \mu }}=\{\mu_1,\ldots,\mu_k\}$. Then, $\widetilde{\alpha_0}$ is a choice set for $\{ \beta_e\mid e\in E\}$ with $\widetilde{\alpha_0}\cap \beta_{(j,q+1)}= \{\mu_j\}$ for all $j\in \{1,\ldots, k\}$ and $\widetilde{\alpha_0}\cap \beta_e= \{\nu_e\}$ for $e\in E\setminus E_{q+1}$. Consequently, $\widetilde{\alpha_0}$ is finite and, thus, $\alpha_0$ belongs to $J_1$. To be more precise, we can even say that ${\rm sym}_\G(\alpha_0)\supseteq \G(P\cup P_{\mbmss{ \mu}})$ holds. $\alpha_0$ is here  the predicate $\alpha^{I,P}_{\mbmss{ \mu} }$  also considered in Table \ref{overview_nota_2}. It is similar to the corresponding predicate in \cite{Gass24B}, but it is not the same. 
 Let $B_{(1,1)}, \ldots, B_{(k,q+1)}$ be new variables for $1$-ary predicates that do not occur in $H( x ,D)$. For any $f\in {\rm assgn}(\S)$, let $\bar f=f\langle {B_{(1,1)} \atop\beta_{(1,1)} } {\cdots\atop\cdots} {B_{(k,q+1)}\atop\beta_{(k,q+1)}}\rangle $. Then, for any $\xi \in I $, any $\delta \in J_1$, and each $\pi\in\G(P)$, we have 
\[\!\!\!\! \!\!\S_{\bar f\langle{x \atop \xi}{D\atop\delta}\rangle} (B_ex \land H( x ,D))=\S_{\bar f\langle{x \atop \pi(\xi)}{D\atop\delta^\pi}\rangle} (B_e x\land H(x ,D))\] 
for every $e\in E $. For any $e\in E$, let $\xi_e$ be the only individual in $\widetilde{\alpha_0}\cap \beta_e$. For $j\in \{1,\ldots,k\}$, let $\alpha_j=I_j$ and $\alpha_{j,0}=\alpha_0\cap \alpha_j$.

By assumption and according to a proof of Lemma \ref{endlSigma}, there is again a choice function $\phi_s$ on $({\cal A}^H_{\xi})_{\xi \in \alpha_0}$ for $s=|\alpha_0|-1$ recursively definable such that the restriction of $\phi_s$ to $\alpha_{j,0}$ denoted by $(\phi_s)|_{\alpha_{j,0}}$ can be the function $\phi_{j,s_j}$ considered in the description of Idea \ref{Idee2}. For each $\xi \in \alpha_0$, let the predicate $\delta_{\xi}$ be defined by $\delta_{\xi}=\phi_s( {\cal A}^H_{\xi})$. Then, for all $ e \in E $, $\delta_{\xi_e} $ is according to the definition of $ {\cal A}^H_{ \xi_e}$ in $\S$ and thus ${\rm sym}_\G(\delta_{ \xi_e})\supseteq\G(P_{\delta_{ \xi_e}})$ holds for some finite $P_{\delta_{ \xi_e}}\subseteq I$. For $\sigma_e^{(0)}=_{\rm df}\{( \xi_e,\eta )\in \bbbn^{1+1}\mid \eta \in \delta_{ \xi_e}\}$, we have ${\rm sym}_{\G}(\sigma_e^{(0)})\supseteq \G(P\cup P_{\mbmss{ \mu}}\cup P_{\delta_{ \xi_e}})$ which means again ${\rm sym}_{\G}(\sigma_0)\supseteq\G(P\cup P_{\mbmss{ \mu}}\cup P_0 )$ for 
$\sigma_0=_{\rm df} \bigcup_{e \in E}\, \sigma_e^{(0)}$ and $P_0=_{\rm df}\bigcup_{e \in E}\, P_{\delta_{ \xi_e}}$.
We want to use
 \[\S\models_{f\langle{A\atop\alpha_0}{S\atop \sigma_0} \rangle}\forall x \exists D (A x \to \forall y (D y \leftrightarrow S xy) \land H(x, D))\] 
and we want to define the permutation $\pi_{\xi}$ as generally as possible and in accordance with the formula $swap_e(x_0, y_0, x, y)$ defined as above for $e \in E $ (here given by $e= (j,i)$ for some $j$ and some $i$). We will use the formula only if $x$ and $x_0$ stand for individuals in $\beta_e$. Consequently, in the combination with $B_ex$, $swap_e(x_0, y_0, x, y)$ will only used for describing permutations in $\G$. $\sigma_e =_{\rm df} \{(\xi , \eta)\in I_j \times I \mid \eta \in (\delta_{ \xi_e})^{\pi_{\xi}}\}$ is again definable by the formula $G_e(x, y)$ defined as above
such that here ${\rm sym}_{\G}(\sigma_e )\supseteq \G(P\cup P_{\mbmss{ \mu}}\cup P_{\delta_{\xi_e} })$ holds. 
Consequently, for $\sigma= \bigcup_{e \in E }\sigma_e$, ${\rm sym}_{\G}(\sigma)\supseteq \G(P\cup P{\mbmss{\mu}}\cup P_0)$ is satisfied and moreover 
\[\S\models_{f\langle{S\atop \sigma} \rangle}\forall x \exists D ( \forall y (D y \leftrightarrow S x y ) \land H( x , D ))\] 
is valid. This means that the consequent of $choice_h^{1,1}(H)$ follows from the antecedent of $choice_h^{1,1}(H)$. 
\hfill \qed

\subsection{Consequences: The independence of TR from AC}\label{sectionConseq}

\begin{corollary} $k\S_0 \models {\rm HAC}$ holds for all $k\geq 1$.
\end{corollary}

\begin{proposition} $k\S_0\models \neg TR^n$ holds for all $k\geq 2$ and for all $n\geq 1$.
\end{proposition}
 
\noindent \noindent {\bf Proof.} Let $k\geq 2$ and $n=1$.  We want to show that there is no binary predicate that belongs to $k\S$ and that is an injective mapping of the infinite set $I_1$ onto a subset of $I_2$ or an injective mapping of $I_2$ into $I_1$.
Without loss of generality, let us assume that $\tau$ is an injective mapping of $I_1$ into  $I_2$ and  that ${\rm sym}_\G(\tau)\supseteq \G(P)$ holds for some finite $P\subseteq I$. Let   $f:I_1\to I_2$ be the  injective mapping  such that  $\eta=f(\xi)$ holds if and only if  $(\xi,\eta)\in \tau$ holds. We use that, for $f(I_1)=_{\rm df}\{f(\xi) \in I_2\mid \xi\in I_1\}$, the set $f(I_1)\setminus P$  is not empty (because we assume that $\bbbn$ and thus $I_1$ are infinite).  Then, there are an $\eta_0\in f(I_1) \setminus  P$, a $\xi_0=f^{-1}(\eta_0)$,  an $\eta \in I_2 \setminus (P\cup \{\eta_0\})$, and a $\pi\in \G(P)$ such that $\xi_0=\pi(\xi_0)$ and  $\eta=\pi(\eta_0)$ hold which implies $f(\xi_0)=\eta_0$ and $f(\xi_0)=\eta$. This contradicts the fact that $f$ and $\tau$ should be  an injective mapping.
 Moreover, for the sets $\{(1,1)\}^{n-1}\times (\{1\}\times \bbbn) $ and $\{(1,1)\}^{n-1}\times (\{2\}\times \bbbn) $, there are also no predicates that could be suitable injective mappings.
\hfill\qed

\vspace{0.3cm}

More general, we have the following.
\begin{proposition} For all $n\geq 1$, $TR^n$ is {\rm HPL}-independent from $choice_h^{1,1}$.
\end{proposition}

By using the ideas for proving Proposition 4.14 in \cite{Gass24B} as described above, it is possible to generalize Proposition \ref{ChoiceInS2}. In this way,  the {\rm HPL}-independence of $TR^r$  from all Ackermann axioms $choice_h^{n,m}(H)$ can be shown for each $r\geq 1$, any $n,m\geq 1$, and any $H\in {\cal L}^{(2)}_{\mbmss{x},D}$.

\begin{theorem}For all $n\geq 1$, $TR^n$ is {\rm HPL}-independent from $choice_h^{(2)}$.
\end{theorem}

As mentioned above, it is not necessary to use ZFC as metatheory and to assume its consistence (for more details see Remarks \ref{WOinFinStr} and \ref{WOinInfinStr}).

\section{An open problem and a remark}\label{section5}

The form of Zorn's lemma $ZL$ considered here states that every set (given by an $n$-ary predicate), that is partially ordered by a relation (given by a $2n$-ary predicate) $\tau$ and whose subsets that are linearly ordered by $\tau$ have an upper bound in this set, has a maximal element with respect to $\tau$. By \cite{Jech73}, this formulation goes also back to Max Zorn (1935). Zorn stated without proof that $ZL$ implies the axiom of choice (for more details see e.g.\,\,\cite[p.\,11]{6}). 

\vspace{0.2cm}
\noindent{\bf Zorn's lemma for $n$-ary predicates}

\nopagebreak 
\noindent\fbox{\parbox{11.8cm}{

\vspace*{0.1cm}

$\begin{array}{ll}
ZL^{n} =_{\rm df} & \forall A \forall T\big(po(T,A) \\ &\qquad\land \forall B ([B\subseteq A]\land \forall \mbm{x}_{1} \forall \mbm{x}_{2} (B\mbm{x}_{1} \land B\mbm{x}_{2} \to T\mbm{x}_{1} \mbm{x}_{2}  \lor T\mbm{x}_{2} \mbm{x}_{1})\\ 
& \qquad\qquad\to \exists \mbm{x}(sup(\mbm{x},B,T,A)))\\ 
&\qquad\to \exists \mbm{x}_{0}(A\mbm{x}_{0} \land \forall \mbm{x}_{1} (A \mbm{x}_{1}\land T\mbm{x}_{0} \mbm{x}_{1}\to \mbm{x}_{0} =\mbm{x}_{1}))).
\end{array}$
}}

\vspace*{0.3cm}

\noindent Here, $\mbm{x}$ and $\mbm{x}_1,\mbm{x}_2,\mbm{x}_3$ are $n$-tuples of variables. $[B \subseteq A]$ stands for $\forall \mbm{x} (B\mbm{x} \to A\mbm{x})$, the formula $po(T,A)$ describes the properties of a partial ordering, and $sup(\mbm{x}_0,B,T,A)$ is the formula $A\mbm{x}_0\land \forall \mbm{x}(B\mbm{x} \to T${\boldmath $xx$}$_0)\land \forall \mbm{x}_{1} (A \mbm{x}_{1} \land \forall \mbm{x}(B\mbm{x} \to T\mbm{x}\mbm{x}_1) \to T\mbm{x}_{0}\mbm{x}_{1})$. $T$ could stand, for instance, for a lexicographic order derived from an order on the set of all individuals.

By \cite[p.\,128]{1c} where $ZL^1$ is denoted by $ZL^{(2)}$, Zorn's lemma $ZL^1$ can be derived from $AC^{1,1}$ in second-order logic with Henkin interpretation. We also know that $AC^{1,1}$ and $TR^1$ are HPL-independent from $ZL^1$ (cf. \cite{Gass84, Gass94}). Whereas we have answered a question posed by Asser for $ZL^1$ by showing that $ZL^1\to TR^1$ does not hold in {\rm HPL}, we do not know whether $TR^1\to ZL^1$ holds in HPL.

\begin{question}
\[^h ax^{(2)}\vdash TR^1\to ZL^1 \mbox{?}\]
\end{question}

Let us remark that, whereas $TR^n \to TR^1$ and $ZL^n \to ZL^1$ are true for all $n\geq 1$ in HPL, both implications $TR^1 \to TR^n$ and $ZL^1 \to ZL^n$ are false for all $n\geq 2$ in HPL (by \cite[p.\,43, Section IV, Satz 2.3.3]{Gass84} and \cite[p.\,539]{Gass94}).

\vspace{0.5cm}
\noindent {\Large\bf Acknowledgment}
\addcontentsline{toc}{section}{\bf Acknowledgment}

I thanks Michael Rathjen and Stephen Mackereth for their questions about why HAC holds in the models analyzed in this paper. The questions were helpful for me to understand how important it is to explain the proof ideas in more detail. 

\end{document}